\documentclass{article}
\usepackage[utf8]{inputenc} 
\usepackage{graphicx} 
\usepackage{booktabs} 
\usepackage{array} 
\usepackage{arydshln}
\usepackage{xcolor}
\usepackage{enumerate}
\usepackage{amsmath}
\usepackage{mathtools}
\usepackage{amssymb}
\usepackage{graphicx}
\usepackage{float}
\usepackage{subfig}
\usepackage{fancyhdr}
\usepackage{dcolumn}
\usepackage[hidelinks]{hyperref}
\usepackage{tikz}
\usetikzlibrary{arrows}
\usetikzlibrary{matrix}
\usetikzlibrary{decorations.pathreplacing}
\usetikzlibrary{calc}
\newcommand{\KW}{\noindent{\bf Keywords: }}
\newcommand{\AMS}{\medskip\noindent{\bf Mathematics Subject Classification 2010: }}

\makeatletter
\newcommand*\dashline{\rotatebox[origin=c]{90}{$\dabar@\dabar@\dabar@$}}
\makeatother

\makeatletter
\newlength\xvec@height%
\newlength\xvec@depth%
\newlength\xvec@width%
\newcommand{\xvec}[2][]{%
  \ifmmode%
    \settoheight{\xvec@height}{$#2$}%
    \settodepth{\xvec@depth}{$#2$}%
    \settowidth{\xvec@width}{$#2$}%
  \else%
    \settoheight{\xvec@height}{#2}%
    \settodepth{\xvec@depth}{#2}%
    \settowidth{\xvec@width}{#2}%
  \fi%
  \def\xvec@arg{#1}%
  \def\xvec@dd{:}%
  \def\xvec@d{.}%
  \raisebox{.2ex}{\raisebox{\xvec@height}{\rlap{%
    \kern.05em
    \begin{tikzpicture}[scale=1]
    \pgfsetroundcap
    \draw (.05em,0)--(\xvec@width-.05em,0);
    \draw (\xvec@width/2,-0.15em)--(\xvec@width/2, 0.15em);
    \ifx\xvec@arg\xvec@d%
      \fill(\xvec@width*.45,.5ex) circle (.5pt);%
    \else\ifx\xvec@arg\xvec@dd%
      \fill(\xvec@width*.30,.5ex) circle (.5pt);%
      \fill(\xvec@width*.65,.5ex) circle (.5pt);%
    \fi\fi%
    \end{tikzpicture}%
  }}}%
  #2%
}
\newcommand{\xvecsub}[2][]{%
  \ifmmode%
    \settoheight{\xvec@height}{$#2$}%
    \settodepth{\xvec@depth}{$#2$}%
    \settowidth{\xvec@width}{$#2$}%
  \else%
    \settoheight{\xvec@height}{#2}%
    \settodepth{\xvec@depth}{#2}%
    \settowidth{\xvec@width}{#2}%
  \fi%
  \def\xvec@arg{#1}%
  \def\xvec@dd{:}%
  \def\xvec@d{.}%
  \raisebox{-0.3ex}{\raisebox{\xvec@height}{\rlap{%
    \kern.05em
    \begin{tikzpicture}[scale=1]
    \pgfsetroundcap
    \draw (.05em,0)--(\xvec@width-.05em,0);
    \draw (\xvec@width/2,-0.1em)--(\xvec@width/2, 0.1em);
    \ifx\xvec@arg\xvec@d%
      \fill(\xvec@width*.45,.5ex) circle (.5pt);%
    \else\ifx\xvec@arg\xvec@dd%
      \fill(\xvec@width*.30,.5ex) circle (.5pt);%
      \fill(\xvec@width*.65,.5ex) circle (.5pt);%
    \fi\fi%
    \end{tikzpicture}%
  }}}%
  #2%
}
\makeatother
 
 \DeclareRobustCommand{\naturalto}{%
  \mathrel{
  \ooalign{$\longrightarrow$\cr%
  \kern1.45ex\raise.275ex\hbox{\scalebox{1}[0.522]{$\mid$}}\cr}
  }%
}

 \DeclareRobustCommand{\naturaltolr}{%
  \mathrel{
  \ooalign{$\longleftrightarrow$\cr%
  \kern1.75ex\raise.275ex\hbox{\scalebox{1}[0.522]{$\mid$}}\cr}
  }%
}

\begin{document}
\title{A Canonical Fuzzy Logic} \author{
Osvaldo Skliar\thanks{Universidad Nacional, Costa Rica. E-mail: osvaldoskliar@gmail.com. https://orcid.org/0000-0002-8521-3858.}
\and Sherry Gapper\thanks{Universidad Nacional, Costa Rica. E-mail: sherry.gapper.morrow@una.ac.cr.  https://orcid.org/0000-0003-4920-6977.}
\and Ricardo E. Monge\thanks{Universidad CENFOTEC, Costa Rica. E-mail: rmonge@ucenfotec.ac.cr. https://orcid.org/0000-0002-4321-5410.}}


\maketitle
\begin{abstract}
A presentation is provided of the basic notions and operations of a) the propositional calculus of a variant of fuzzy logic -- canonical fuzzy logic, CFL -- and in a more succinct and introductory way, of b) the theory of fuzzy sets according to that same logic. The propositional calculus of bivalent classical logic and classical set theory can be considered as particular cases of the corresponding theories of CFL if the numerical value of a specific parameter $w$ is restricted to only two possibilities, 0 and 1.
\end{abstract}

\KW  canonical fuzzy logic, classical logic, propositional calculus, set theory, the weight or degree of truth of a proposition, weight or degree of membership of an element in a set

\AMS 03B05, 03B10, 03B52, 03E72

\newpage
\section{Introduction}
The objective of this article is to present the main ideas of a \textit{c}anonical \textit{f}uzzy \textit{l}ogic (CFL).  This logic is a particular case of a non-classical logical formalism -- in the sense that it differs from that of \textit{b}ivalent classical \textit{l}ogic (BL) -- that is susceptible to diverse interpretations each of which is a non-classical logic. In this first treatment of the topic, only one of these interpretations will be considered: precisely that of CFL.

In this article BL will be understood as the mathematical logic in which there are only two possible truth values, true and false, as discussed in outstanding contributions by G. Boole, F. L. G. Frege, C. S. Pierce and B. Russell. This logic can be considered a notable improvement and an extraordinary amplification of Aristotelian logic.

CFL is a parametric type of logic inasmuch as a parameter is used: $w$ -- weight or degree of truth. It is known as canonical because there is a marked continuity with BL: the “laws” -- theorems -- of this logic conserve their validity in CFL. BL can be regarded as a particular ``limit'' case of CFL.

For a clear understanding of this article, knowledge of only basic notions of propositional calculus of BL and of classical set theory is required. For an introductory overview of these topics, one may consult, for example: \cite{a1}, \cite{a2}, \cite{a3} and \cite{a4}. On fuzzy logic and fuzzy set theory, one may consult, for example: \cite{b1}, \cite{b2}, \cite{c1} and \cite{c2}. An interesting view of the development of fuzzy logic and the theory of fuzzy sets, as well as of the prospects of those fields was presented by Lotfi A. Zadeh \cite{c3}.

\section{Some Correspondences between Operations of Propositional Calculus of BL and Operations of Classical Set Theory}
Variables which can be replaced by propositions are known as propositional variables. Letters such as $p$, $q$, $r$, $\ldots$ are commonly used to refer to different propositional variables. In future work on the non-classical logical formalism mentioned above, the letter $p$ will be used as an abbreviation for \textit{p}robability. Thus, to prevent confusion, in this article reference will be made to the different propositional variables as follows: $q_1$, $q_2$, $q_3$, $\ldots$.

Strictly speaking, attributing a truth value to a propositional variable is not admissible. It only makes sense to attribute a truth value -- either true or false – to a proposition. However, in line with a license commonly found in specialized literature, expressions like ``Consider that $q_1$ has been replaced by a true proposition and $q_2$ has been replaced by a false proposition'' will be communicated in an abbreviated way such as ``$q_1$ is true and $q_2$ is false''.

In logical expressions in which only the symbol corresponding to a sole proposition (and eventually its negation) is present, the subscript can be eliminated and the proposition can be symbolized simply as $q$.

Different sets will be symbolized as $C_1$, $C_2$, $C_3$, $\ldots$. In expressions in which reference is made to only one of those sets, and eventually to its complement, the subscript can be eliminated and the set can be  symbolized simply as $C$.

Two very theoretically important sets, the universal set and its complement – the empty set – will be symbolized respectively as follows: $\mathbb U$ and $\varnothing$.

With no exceptions, all the elements to be considered when using set theory for a specific topic belong to the universal set $\mathbb U$, also called universe of discourse or simply universe. Any element belonging to the universal set $\mathbb U$ does not belong to the corresponding empty set  $\varnothing$. Therefore, in classical set theory, no element belongs to $\varnothing$.

Some of the symbols used in set theory to express different relationships between sets are those of equality ($=$), inclusion in a broad sense ($\subseteq$), and inclusion in a narrow sense ($\subset$). 

$C_1 = C_2$ is read as: ``The set $C_1$ is equal to the set $C_2$''. The meaning of $C_1 = C_2$ is as follows: The same elements of a given universal set $\mathbb U$ that belong to $C_1$ belong to $C_2$. In other words, if any element whatsoever of a given universal set $\mathbb U$ belongs to  $C_1$, then it also belongs to  $C_2$, and if any element whatsoever of a given universal set $\mathbb U$ does not belong to  $C_1$, then it does not belong to  $C_2$ either.

$C_1 \subseteq C_2$ is read as: ``$C_1$ is included, in a broad sense, in $C_2$''. The meaning of $C_1 \subseteq C_2$ is as follows: Any element of a given universal set $\mathbb U$ which belongs to $C_1$ also belongs to $C_2$, and the possibility that $C_1$ may be equal to $C_2$ ($C_1 = C_2$) is not excluded.

$C_1 \subset C_2$ is read as: ``$C_1$ is included, in a narrow sense, in $C_2$''. The meaning of $C_1 \subset C_2$ is as follows: Any element of a given universal set $\mathbb U$ which belongs to $C_1$ also belongs to $C_2$, and there is at least one element in that universal set $\mathbb U$ which belongs to  $C_2$ but does not belong to $C_1$. Thus the possibility that $C_1$ is equal to $C_2$ is excluded.

The logical connective of \textit{negation} of a proposition is usually symbolized as $\neg$. What many authors express with $\neg q$ (not $q$) will be expressed in this article by a horizontal bar above the negated proposition as $\overline{q}$.

The operator of complementation of a set $C$ is usually represented by a bar above the symbol $C$: $\overline{C}$. In this article a slightly different symbol will be used to represent that complement. In this way, confusion between the logical connective of negation in propositional calculus and the operator of complementation used in set theory can be prevented. A proposition -- that is, a statement also known as a ``sentence'' at times -- may be negated. However, according to the usual notion of a set, it would be senseless to negate a set. Therefore, the set  $\xvec{C}$ will be called the complement set of $C$. Recall that all the elements belonging to the universal set $\mathbb U$ considered that do not belong to $C$ belong to $\xvec{C}$.

The truth table (a) corresponding to the negation of q -- that is, $\overline{q}$ -- and the membership table (b) corresponding to the complementary set $\xvec{C}$  of any set $C$ are shown in figure \ref{f1}. 

 \begin{figure}[H]
\centering
\subfloat[Truth table corresponding to the negation of q; that is, $\overline{q}$.]{
\hspace{2in}
\begin{tabular}{c||c}
$q$ & $\overline{q}$  \\
\midrule
0 & 1   \\ 
1 & 0   \\ 
\end{tabular}
\hspace{2in}
\label{f1a}
}
\\
\subfloat[Membership table corresponding to the complementary set of $C$ -- that is, $\xvec{C}$. $C$ is any set such that all the elements belonging to it are elements of the universal set $\mathbb U$ considered above.]{
\hspace{1.9in}
\begin{tabular}{c||c}
$C$ & $\xvec{C}$  \\
\midrule
0 & 1   \\ 
1 & 0   \\ 
\end{tabular}
\hspace{2in}
\label{f1b}}%
\caption{a) Truth table corresponding to $\overline{q}$, and b) membership table corresponding to $\xvec{C}$.}
\label{f1}
\end{figure}
 
In a truth table, a zero (0) in the column corresponding to a proposition -- that is, a truth value equal to 0 -- represents the supposition that the proposition is false. A one (1) in the column corresponding to a proposition -- that is, a truth value equal to 1 -- represents the supposition that the proposition is true. 

In the membership table, a zero (0) in the column corresponding to a set -- that is, a membership value equal to 0 -- represents the supposition that any element of the universal set  $\mathbb U$ considered does not belong to that set. A one (1) in the column corresponding to a set -- that is, a membership value equal to 1 -- represents the supposition that any element whatsoever of $\mathbb U$ does belong to that set.

The first row of the truth table in figure \ref{f1a}, in which a sequence of digits (0 and 1) appears, should be interpreted as follows: If $q$ is false, then $\overline{q}$ is true. The second row of that truth table in which a sequence of digits (1 and 0) appears, should be interpreted as follows: If $q$ is true, then $\overline{q}$ is false.

The first row of the membership table in figure \ref{f1b}, in which a sequence of digits (0 and 1) appears, should be interpreted as follows: If  any element belonging to the $\mathbb U$ considered does not belong to $C$, then it does belong to $\xvec{C}$. The second row of the membership table in figure \ref{f1b}, in which a sequence of digits (1 and 0) appears, should be interpreted as follows: If any element belonging to the $\mathbb U$ considered belongs to $C$, then that element does not belong to $\xvec{C}$.

Note the equality of the first and second rows, respectively, of the tables in figures \ref{f1a} and \ref{f1b}. This makes it possible to establish: 1) a correspondence between the logical connective of negation of propositional calculus and the operator of complementation of set theory; 2) a correspondence between the proposition $q$ and the set $C$; and 3) a correspondence between $\overline{q}$ and the set $\xvec{C}$.

In propositional calculus, the logical connective of \textit{disjunction}, known as ``or'', is symbolized as $\lor$. The proposition $q_1 \lor q_2$ is read as ``$q_1$ or $q_2$''. In set theory, the operator of the union of sets is symbolized as $\cup$. The set $C_1 \cup C_2$ is the set resulting from the operation of union of the sets $C_1$ and $C_2$. It is admitted that each element belonging to set $C_1$, and each element belonging to set $C_2$ are elements belonging to the universal set $\mathbb U$, which should be specified unambiguously. This is achieved by specifying the elements belonging to $\mathbb U$. 

Figure \ref{f2} presents a) the truth table corresponding to $q_1 \lor q_2$, and the membership table corresponding to $C_1 \cup C_2$.

\begin{figure}[H]
\centering
\subfloat[Truth table for the disjunction of the propositions $q_1$ and $q_2$: $q_1 \lor q_2$.]{
\hspace{1.6in}
\begin{tabular}{c|c||c}
$q_1$ &  $q_2$ & $q_1 \lor q_2$  \\
\midrule
0 & 0 & 0   \\  
0 & 1 & 1   \\  
1 & 0 & 1   \\  
1 & 1 & 1   \\  
\end{tabular}
\hspace{2in}
\label{f2a}
}
\\
\subfloat[Membership table for the union of $C_1$ and $C_2$: $C_1 \cup C_2$.]{
\hspace{1.5in}
\begin{tabular}{c|c||c}
$C_1$ &  $C_2$ & $C_1 \cup C_2$  \\
\midrule
0 & 0 & 0   \\  
0 & 1 & 1   \\  
1 & 0 & 1   \\  
1 & 1 & 1   \\  
\end{tabular}
\hspace{2in}
\label{f2b}}%
\caption{a) Truth table corresponding to $q_1 \lor q_2$; and b) membership table corresponding to $C_1 \cup C_2$.}
\label{f2}
\end{figure}

The sequence of digits 0, 0, 0 in the first rows of the tables in figures \ref{f2a} and \ref{f2b} should be interpreted as follows: In figure \ref{f2a}, that sequence means that if $q_1$ is false, and $q_2$ also is false, then the proposition $q_1 \lor q_2$ is false. In figure \ref{f2b}, that same sequence of digits means that if any element of the $\mathbb U$ considered belongs neither to $C_1$ nor to $C_2$, then that element does not belong to the set $C_1 \cup C_2$.

The sequence of digits 0, 1, 1 in the second rows of the tables in figures \ref{f2a} and \ref{f2b} should be interpreted as follows: In figure \ref{f2a}, that sequence means that if $q_1$ is false, and $q_2$ is true, then the proposition $q_1 \lor q_2$ is true. In figure \ref{f2b}, that same sequence of digits means that if any element of the $\mathbb U$ considered does not belong to $C_1$, and that element does belong to $C_2$, then that element does belong to the set $C_1 \cup C_2$.

The sequence of digits 1, 0, 1 in the third rows of the tables in figures 2a and 2b should be interpreted as follows: In figure \ref{f2a}, that sequence means that if $q_1$ is true and $q_2$ is false, then the proposition $q_1 \lor q_2$ is true. In figure \ref{f2b}, that same sequence of digits means that if any element of the $\mathbb U$ considered belongs to $C_1$, and that element does not belong to $C_2$, then that element does belong to $C_1 \cup C_2$.

The sequence of digits 1, 1, 1 in the fourth rows of the tables in figures \ref{f2a} and \ref{f2b} should be interpreted as follows: In figure \ref{f2a}, that sequence means that if $q_1$ is true, and $q_2$ is also true, then $q_1 \lor q_2$ is true. In figure \ref{f2b}, that same sequence of digits means that if any element of the $\mathbb U$ considered belongs to $C_1$, and that element also belongs to $C_2$, then that element belongs to $C_1 \cup C_2$.

Observe that only when both the proposition $q_1$ and the proposition $q_2$ are false, the proposition $q_1 \lor q_2$ is false. Also, only when any element of the universal set $\mathbb U$ does not belong to $C_1$ and does not belong to $C_2$, then that element does not belong to $C_1 \cup C_2$.

Note the equality of the first, second, third and fourth rows respectively, of the tables in figures \ref{f2a} and \ref{f2b}. This equality makes it possible to establish the following correspondences: 1) the propositions $q_1$ and $q_2$ correspond, respectively, to the sets $C_1$ and $C_2$; and 2) the logical connective of disjunction $\lor$ of propositional calculus corresponds to the union operator $\cup$ in set theory.

Detailed explanations like these for each pair of tables presented in figures \ref{f3}, \ref{f4} and \ref{f5} might be tedious to read. Briefer explanations will be provided for each of those pairs of tables, and based on the above information, the missing details can be completed. Also for each pair of tables, correspondences may be established between: 1) the propositions $q_1$ and $q_2$ and the sets $C_1$ and $C_2$, respectively; and 2) a particular logical connective of propositional calculus and a particular operator of set theory.

In propositional calculus, the logical connective of \textit{conjunction} ``and'' is symbolized as $\land$.  The proposition $q_1 \land q_2$ is read as ``$q_1$ and $q_2$''. In set theory, the operator of intersection is symbolized as $\cap$. The set $C_1 \cap C_2$ is the set resulting from the operation of intersection between the sets $C_1$ and $C_2$. It is admitted that each element belonging to the set $C_1$ and each element belonging to the set $C_2$ are elements belonging to a universal set $\mathbb U$ which should be specified unambiguously. This is achieved by specifying the elements in $\mathbb U$, that is, the elements belonging to the set $\mathbb U$.

The truth table (a) corresponding to $q_1 \land q_2$, and the membership table (b) corresponding to $C_1 \cap C_2$ are shown in figure \ref{f3}.

\begin{figure}[H]
\centering
\subfloat[Truth table for the conjunction of any two propositions $q_1$ and  $q_2$: $q_1 \land q_2$.]{
\hspace{1.6in}
\begin{tabular}{c|c||c}
$q_1$ &  $q_2$ & $q_1 \land q_2$  \\
\midrule
0 & 0 & 0   \\  
0 & 1 & 0   \\  
1 & 0 & 0   \\  
1 & 1 & 1   \\  
\end{tabular}
\hspace{2in}
\label{f3a}
}
\\
\subfloat[Membership table corresponding to the intersection of $C_1$ and $C_2$: $C_1 \cap C_2$.]{
\hspace{1.5in}
\begin{tabular}{c|c||c}
$C_1$ &  $C_2$ & $C_1 \cap C_2$  \\
\midrule
0 & 0 & 0   \\  
0 & 1 & 0   \\  
1 & 0 & 0   \\  
1 & 1 & 1   \\  
\end{tabular}
\hspace{2in}
\label{f3b}}%
\caption{a) Truth table corresponding to $q_1 \land q_2$; and b) membership table corresponding to $C_1 \cap C_2$.}
\label{f3}
\end{figure}

In figure \ref{f3a}, it can be seen that only when $q_1$ is true and  $q_2$ is also true, the conjunction of them $q_1 \land q_2$ is true. This is the case considered in the fourth row of the corresponding truth table. In the other three possible cases, in the other three rows of that truth table, $q_1 \land q_2$ is false. Likewise, in figure \ref{f3b}, it is seen that only when any element of $\mathbb U$ that belongs to $C_1$ and also belongs to $C_2$, that element will belong to $C_1 \cap C_2$. This case is considered in the fourth row of the respective membership table. In the other three cases, in the other three rows of that membership table, that element does not belong to $C_1 \cap C_2$.

The consideration of the tables presented in figure \ref{f3} makes it possible to establish the following correspondences: 1) a correspondence between the propositions $q_1$ and $q_2$ and the sets $C_1$ and  $C_2$ respectively; and 2) a correspondence between the logical connective of conjunction $\land$ of propositional calculus, and the operator of intersection $\cap$ in set theory.

In propositional calculus, the logical connective of \textit{material implication} is symbolized as $\to$.  The proposition $q_1 \to q_2$, considered as conditional, is read as ``if $q_1$, then $q_2$''. The proposition $q_2 \to q_1$ is read as ``if $q_2$, then $q_1$''.  The antecedent of the proposition $q_1 \to q_2$ is $q_1$, and the consequent of that proposition is $q_2$. Likewise,  the antecedent of the proposition $q_2 \to q_1$ is $q_2$, and the consequent of that proposition is $q_1$.

In set theory, the operator of the material implication of sets is symbolized as $\naturalto$. The set $C_1 \naturalto C_2$ is generated by applying the operator of material implication to the ordered pair of sets ($C_1$, $C_2$). The set $C_2 \naturalto C_1$ is generated by applying the operator of material implication to the ordered pair of sets ($C_2$, $C_1$).

Figure \ref{f4} presents (a) the truth tables corresponding to the propositions $q_1 \to q_2$ and $q_2 \to q_1$, and (b) the membership tables corresponding to the sets $C_1 \naturalto C_2$ and $C_2 \naturalto C_1$.

\begin{figure}[H]
\centering
\subfloat[Truth table corresponding to the proposition $q_1 \to q_2$ and truth table corresponding to $q_2 \to q_1$.]{
\hspace{1in}
\begin{tabular}{c|c||c|c}
$q_1$ &  $q_2$ & $q_1 \to q_2$  & $q_2 \to q_1$  \\
\midrule
0 & 0 & 1 & 1  \\  
0 & 1 & 1 & 0   \\  
1 & 0 & 0 & 1 \\  
1 & 1 & 1 & 1   \\  
\end{tabular}
\hspace{1in}
\label{f4a}
}
\\
\subfloat[Membership table corresponding to the set $C_1 \naturalto C_2$ and membership table corresponding to the set $C_2 \naturalto C_1$.]{
\hspace{1in}
\begin{tabular}{c|c||c|c} 
$C_1$ &  $C_2$ & $C_1 \naturalto C_2$  & $C_2 \naturalto C_1$ \\
\midrule
0 & 0 & 1 & 1  \\  
0 & 1 & 1 & 0   \\  
1 & 0 & 0 & 1 \\  
1 & 1 & 1 & 1   \\  
\end{tabular}
\hspace{1in}
\label{f4b}}%
\caption{a) Truth tables corresponding to the propositions $q_1 \to q_2$ and $q_2 \to q_1$, and membership tables corresponding to the sets $C_1 \naturalto C_2$ and $C_2 \naturalto C_1$.}
\label{f4}
\end{figure}

Note that the only case in which the proposition $q_1 \to q_2$ is false is that in which the antecedent $q_1$ is true and the consequent $q_2$ is false. This is the case of the third row in the truth table $q_1 \to q_2$. In the other three possible cases, $q_1 \to q_2$ is true. Likewise, the only case in which the proposition $q_2 \to q_1$ is false is that in which the antecedent $q_2$ is true and the consequent $q_1$ is false. This is the case of the second row of the truth table for the proposition $q_2 \to q_1$. In the other three cases, $q_2 \to q_1$ is true.

In addition, notice also that the only case in which any element of the universal set $\mathbb U$ considered does not belong to the set $C_1 \naturalto C_2$ is that in which that element belongs to $C_1$ and does not belong to  $C_2$. This is the case of the third row of the membership table for $C_1 \naturalto C_2$. Likewise, the only case in which any element of the universal set $\mathbb U$ considered does not belong to the set $C_2 \naturalto C_1$ is that in which that element belongs to $C_2$ and does not belong to  $C_1$.

The consideration given to the tables in figure \ref{f4} makes it possible to establish the following correspondences: 1) a correspondence between the propositions $q_1$ and $q_2$ and the sets $C_1$ and $C_2$, respectively; and 2) a correspondence between the connective of material implication of propositional calculus and the operator of material implication of set theory. 

In propositional calculus, the connective of \textit{material bi-implication} (or logical equivalence) is symbolized as: $\longleftrightarrow$.  The proposition $q_1 \longleftrightarrow q_2$ is true if and only if $q_1$ and $q_2$ have the same truth value; that is, if both are true or both are false. In set theory, the operator of material bi-implication is symbolized as $\naturaltolr$. Any element of the universal set $\mathbb U$ that belongs both to $C_1$ and $C_2$ belongs to $C_1 \naturaltolr C_2$; so does any element that does not belong to either $C_1$ or $C_2$. However, if any element of $\mathbb U$ belongs to only one of the two sets, $C_1$ or $C_2$, and not to the other, then it does not belong to the set $C_1 \naturaltolr C_2$.

The truth table (a) corresponding to the proposition $q_1 \longleftrightarrow q_2$, and the membership table (b) of the set $C_1 \naturaltolr C_2$ are represented in figure \ref{f5}.

\begin{figure}[H]
\centering
\subfloat[Truth table corresponding to the proposition $q_1 \longleftrightarrow q_2$.]{
\hspace{1.5in}
\begin{tabular}{c|c||c}
$q_1$ &  $q_2$ & $q_1 \longleftrightarrow q_2$  \\
\midrule
0 & 0 & 1  \\  
0 & 1 & 0   \\  
1 & 0 & 0   \\  
1 & 1 & 1   \\  
\end{tabular}
\hspace{2in}
\label{f5a}
}
\\
\subfloat[Membership table corresponding to the set $C_1 \naturaltolr C_2$.]{
\hspace{1.4in}
\begin{tabular}{c|c||c}
$C_1$ &  $C_2$ & $C_1 \naturaltolr C_2$  \\
\midrule
0 & 0 & 1  \\  
0 & 1 & 0   \\  
1 & 0 & 0   \\  
1 & 1 & 1   \\  
\end{tabular}
\hspace{2in}
\label{f5b}}%
\caption{a) Truth table corresponding to the proposition $q_1 \longleftrightarrow q_2$; and b) membership table corresponding to the set $C_1 \naturaltolr C_2$.}
\label{f5}
\end{figure}

The comparison of the tables presented in figures \ref{f5a} and  \ref{f5b} makes it possible to establish the following correspondences: 1) a correspondence between the propositions $q_1$ and $q_2$ and the sets $C_1$ and $C_2$, respectively; and 2) a correspondence between the logical connective of material bi-implication of propositional calculus and the operator of material bi-implication of set theory.

The comparison of each of the truth tables and its corresponding membership table presented in this section makes it possible to establish the following correspondences: 1) between each proposition and a set, and 2) between each logical connective of propositional calculus and a specific operator of set theory.

Recall the number of existing logical functions, which are also propositions, of $n$ propositions, for $n=1,2,3\ldots$. If in the logical operations carried out there are $n$ propositions -- $q_1$, $q_2$, $\ldots, q_n$ -- in the corresponding truth tables, there will be $2^n$  rows because each of those propositions can have two truth values -- true or false. Each row of those tables will correspond to each possible case of different assignments for the truth values of each of those $n$ propositions. Since for each of these cases there are two possible assignments of truth value for the logical function to be specified, there are $2^{(2^n)}$ possible logical functions of $n$ propositions. Thus, for example, if $n=2$, there are 16 possible logical functions, and if $n=3$, there are 256 possible logical functions. 

For each of the $2^{(2^n)}$ logical functions of $n$ propositions, for $n=1,2,3\ldots$, there is, according to the approach used, a function of $n$ sets, which also is a set.

For propositions resulting from the use of connectives, such as $q_1 \lor q_2$ or $q_3 \land q_4$, it can be suitable to express them in parentheses as $(q_1 \lor q_2)$ and $(q_3 \land q_4)$, respectively. Thus, if a connective from propositional calculus is used with those propositions to obtain another proposition, it is clear how that has operated. For example, $(q_1 \lor q_2)\to(q_3 \land q_4)$ is the proposition of a conditional nature: ``Si $q_1 \lor q_2$, then $q_3 \land q_4$'', in which $q_1 \lor q_2$ is the antecedent, and $q_3 \land q_4$ is the consequent. The proposition  $(q_1 \lor q_2)\to(q_3 \land q_4)$ has been obtained by the action of the connective of material implication on the following ordered pair of propositions $\left( (q_1 \lor q_2),(q_3 \land q_4)\right)$. Likewise, for clarity, the proposition $(q_1 \lor q_2)\to(q_3 \land q_4)$ can be expressed between parentheses if an operation is carried out on it and on some other proposition, by using some logical connective. Hence, for example, $\left( (q_1 \lor q_2)\to(q_3 \land q_4) \right) \lor (q_4 \to q_5)$ is the proposition obtained through the action of the logical connective of disjunction on the propositions $(q_1 \lor q_2)\to(q_3 \land q_4)$, and $(q_4 \to q_5)$.

Given the correspondences mentioned, 1) between propositions and sets, and 2) between logical connectives and set theory operators, considerations of this same type concerning the use of parentheses are valid in this theory. Therefore, the set $(C_1 \cup C_2) \naturalto (C_3 \cap C_4)$ corresponds to $(q_1 \lor q_2)\to(q_3 \land q_4)$; and the set $\left((C_1 \cup C_2) \naturalto (C_3 \cap C_4) \right) \cup (C_4 \naturalto C_5)$ corresponds to $\left( (q_1 \lor q_2)\to(q_3 \land q_4) \right) \lor (q_4 \to q_5)$.

\section{Isomorphism between Each Law or Theorem of Propositional Calculus and the Corresponding Expression of the Universal Set}

If a function of $n$ propositions, for $n = 1, 2, 3, \ldots$, is true, regardless of the truth values of each of those $n$ propositions, then that propositional function, which also is a proposition, is considered a law -- a theorem -- of propositional calculus. 

When carrying out operations with $n$ sets, for $n = 1, 2, 3, \ldots$, using the operators mentioned above, it will be supposed that each element belonging to each of those sets is an element belonging to the universal set $\mathbb U$ specified unambiguously. If the set resulting from these operations is equal to the universal set $\mathbb U$, regardless of what those sets are (with the only limitation being that specified above), then the equation which expresses the equality between the resulting set mentioned and the universal set $\mathbb U$ is considered a law -- a theorem -- of set theory.

Some examples of laws of propositional calculus and the corresponding laws of set theory will be considered below.

The law (or the principle) of the excluded middle, known also as the law (or the principle) of the excluded third (in Latin, \textit{principium tertii exclusi}) can be expressed as $q \lor \overline{q}$. Another Latin designation for that law is \textit{tertium non datur}: no third (possibility) is given.

The expression of a set which is isomorphic to that law is $C \cup {\xvec{C}}$.

In figure \ref{f6}, the truth table (a) corresponding to the law of the excluded middle, and the membership table (b) of the set whose expression is isomorphic to that law of propositional calculus are presented.

 \begin{figure}[H]
\centering
\subfloat[Truth table corresponding to the law of propositional calculus $q\lor \overline{q} $.]{
\hspace{1.7in}
\begin{tabular}{c|c||c}
$q$ & $\overline{q}$  &  $q\lor \overline{q} $\\
\midrule
0 & 1 & 1   \\  
1 & 0  & 1 \\  
\end{tabular}
\hspace{2in}
\label{f6a}
}
\\
\subfloat[Membership table corresponding to the set $C \cup {\xvec{C}}$.]{
\hspace{1.6in}
\begin{tabular}{c|c||c}
$C$ & $\xvec{C}$  & $C \cup \xvec{C}$ \\
\midrule
0 & 1   & 1\\  
1 & 0   & 1 \\  
\end{tabular}
\hspace{2in}
\label{f6b}}%
\caption{a) Truth table for $q\lor \overline{q}$, and b) membership table for $C \cup {\xvec{C}}$.}
\label{f6}
\end{figure}

In the truth table in figure \ref{f6a}, it can be seen that both if $q$ is false and if $q$ is true, $q\lor \overline{q}$ is true. In other words, $q\lor \overline{q}$ is true due to its logical form, regardless of whether $q$ is false or $q$ is true.

In the first row of the membership table for the set  $C \cup {\xvec{C}}$, there is a sequence of digits 0, 1, 1. This sequence should be interpreted as follows: If any element of $\mathbb{U}$ does not belong to $C$, then, given the definition of the complement of $C$ (${\xvec{C}}$), that element belongs to ${\xvec{C}}$, and therefore, given the characterization of $C \cup {\xvec{C}}$, that element belongs to  $C \cup {\xvec{C}}$.

In the second row of the membership table for the set $C \cup {\xvec{C}}$, there is a sequence of digits 1, 0, 1. This sequence should be interpreted as follows: If any element of $\mathbb{U}$ belongs to $C$, then, given the definition of the complement of $C$ (${\xvec{C}}$), that element does not belong to ${\xvec{C}}$, and therefore, given the characterization of $C \cup {\xvec{C}}$, that element belongs to  $C \cup {\xvec{C}}$.

It can be seen, then, that any element of $\mathbb{U}$ belongs to $C \cup {\xvec{C}}$, regardless of whether it belongs to $C$, or does not belong to $C$. Therefore, $C \cup {\xvec{C}}$ is an expression of the set $\mathbb{U}$: $C \cup {\xvec{C}} = \mathbb{U}$.

Notice that $C \cup {\xvec{C}}$ is an expression of $\mathbb{U}$, which is isomorphic to  $q\lor \overline{q}$. Indeed, $C$ corresponds to $q$, ${\xvec{C}}$ corresponds to $\overline{q}$ and the operator for set union  -- $\cup$ -- corresponds to the connective $\lor$. The equality $C \cup {\xvec{C}}$ = $\mathbb{U}$ is a law, or theorem, of set theory.

The law known as \textit{modus tolendo tolens} of propositional calculus can be expressed as $\left( (q_1 \to q_2) \land \overline{q}_2 \right) \to q_1$. The expression of a set which is isomorphic to that law of propositional calculus is: $( (C_1 \naturalto C_2) \cap {\xvec{C}_2} ) \naturalto {\xvec{C}_1}$.

The truth table (a) corresponding to the law \textit{modus tolendo tolens} of propositional calculus and the membership table (b) for a set whose expression is isomorphic to that law are presented in figure \ref{f7}. 

\begin{figure}[H]
\centering
\subfloat[Truth table for the law of propositional calculus $((q_1 \to q_2) \land  \overline{q}_2) \to \overline{q}_1 $.]{
\begin{tabular}{c|c|c|c|c|c||c} 
$q_1$ &  $q_2$ & $\overline{q}_1$  & $\overline{q}_2$  & $q_1 \to q_2$  & $ (q_1 \to q_2) \land \overline{q}_2 $ & $((q_1 \to q_2) \land \overline{q}_2) \to \overline{q}_1$\\
\midrule
0 & 0 & 1 & 1 & 1 & 1 & 1 \\  
0 & 1 & 1 & 0 & 1 & 0 & 1   \\  
1 & 0 & 0 & 1 & 0 & 0 & 1   \\  
1 & 1 & 0 & 0 & 1 & 0 & 1   \\  
\end{tabular}
\label{f7a}
}
\\
\subfloat[Membership table for the expression of a set -- $ ((C_1 \naturalto C_2) \cap {\xvec{C}_2} ) \naturalto {\xvec{C}_1}$ -- isomorphic to the law of propositional calculus $((q_1 \to q_2) \land \overline{q}_2)\to \overline{q}_1$.]{
\hspace{-.3in}
\begin{tabular}{c|c|c|c|c|c||c}
$C_1$ &  $C_2$ & ${\xvec{C}_1}$ & ${\xvec{C}_2}$ & $C_1 \naturalto C_2$ & $ (C_1 \naturalto C_2) \cap {\xvec{C}_2} $  & $( (C_1 \naturalto C_2) \cap {\xvec{C}_2} ) \naturalto {\xvec{C}_1} $\\
\midrule
0 & 0 & 1 & 1 & 1 & 1 & 1 \\  
0 & 1 & 1 & 0 & 1 & 0 & 1   \\  
1 & 0 & 0 & 1 & 0 & 0 & 1   \\  
1 & 1 & 0 & 0 & 1 & 0 & 1   \\  
\end{tabular}
\label{f7b}%
}
\caption{a) Truth table corresponding to the proposition $((q_1 \to q_2) \land \overline{q}_2)\to \overline{q}_1$, and b) membership table corresponding to $( (C_1 \naturalto C_2) \cap {\xvec{C}_2} ) \naturalto {\xvec{C}_1} $.}
\label{f7}
\end{figure}

It can be seen in figure \ref{f7} that the truth table and the membership table considered were represented to the right of the two parallel segments of vertical lines. For obvious reasons, knowledge of the tables represented on the left of the tables is essential to obtain these tables.

It can be observed that  $((q_1 \to q_2) \land \overline{q}_2)\to \overline{q}_1$ is a tautology of propositional calculus, given that it is true regardless of the truth values of $q_1$ and $q_2$. Moreover, given that any element of $\mathbb U$ belongs to the set $((C_1 \naturalto C_2) \cap \xvec{C}_2 ) \naturalto \xvec{C}_1$, this set is equal to the universal set $(((C_1 \naturalto C_2) \cap \xvec{C}_2 ) \naturalto \xvec{C}_1) = \mathbb U$.

The tautology and the expression of the universal sets $\mathbb U$ considered are isomorphic. Indeed, $q_1$ corresponds to $C_1$,  $\overline{q}_1$ corresponds to ${\xvec{C}_1}$, $q_2$ corresponds to $C_2$,  $\overline{q}_2$ corresponds to ${\xvec{C}_2}$, and the logical connectives of conjunction $\land$ and of material implication $\to$ of propositional calculus correspond, respectively, to the operators of intersection $\cap$ and of the material implication $\naturalto$ of sets.

Reference will be made below to some other laws of propositional calculus and to the corresponding expressions of $\mathbb U$ which are isomorphic to them. Given the information provided above, it is clear how to construct the corresponding truth tables and membership tables. For this reason, those tables are not included here.

The law known as \textit{modus ponendo ponens} can be expressed as $((q_1 \to q_2)\land q_1) \to q_2$. This proposition is true regardless of the truth values of $q_1$ and $q_2$. It is then a tautology. The corresponding expression which is isomorphic to it is $((C_1 \naturalto C_2)\cap C_1) \naturalto C_2$. The corresponding law, or theorem, in set theory is $(((C_1 \naturalto C_2)\cap C_1) \naturalto C_2)=\mathbb U$. 

The law of transitivity of propositional calculus can be expressed as: $((q_1 \to q_2) \land (q_2 \to q_3)) \to (q_1 \to q_3)$. The expression of the universal set which is isomorphic to that law is as follows: $((C_1 \naturalto C_2)\cap (C_2 \naturalto C_3)) \naturalto (C_1 \naturalto C_3)$. The corresponding law in set theory is: $(((C_1 \naturalto C_2)\cap (C_2 \naturalto C_3)) \naturalto (C_1 \naturalto C_3)) = \mathbb{U}$.

One of De Morgan's laws in propositional calculus is: $(\overline{q_1 \lor q_2}) \longleftrightarrow (\overline{q}_1 \land \overline{q}_2)$. The expression of the universal set which is isomorphic to that law is: $({\xvec{C_1 \cup C_2}}) \naturaltolr ({\xvec{C}_1} \cap {\xvec{C}_2})$. The corresponding law in set theory is: $(({\xvec{C_1 \cup C_2}}) \naturaltolr ({\xvec{C}_1} \cap {\xvec{C}_2})) = \mathbb{U}$.

The other of De Morgan's laws in propositional calculus is: $(\overline{q_1 \land q_2}) \longleftrightarrow (\overline{q}_1 \lor \overline{q}_2)$. The expression of the universal set which is isomorphic to that law is: $({\xvec{C_1 \cap C_2}}) \naturaltolr ({\xvec{C}_1} \cup {\xvec{C}_2})$. The corresponding law in set theory is $(({\xvec{C_1 \cap C_2}}) \naturaltolr ({\xvec{C}_1} \cup {\xvec{C}_2})) = \mathbb U$.

If a function of $n$ propositions, for $n = 1, 2, 3, \ldots$ is false, regardless of the truth values of each of these $n$ propositions, then that propositional function, which also is a proposition, is considered a \textit{contradiction}. 

The proposition $q_1 \land \overline q $ is considered below. This proposition, a propositional function of a sole proposition $q$ is a contradiction, according to the above criterion.

The truth table (a) of $q\land \overline{q} $ and the membership table (b) of the expression of a set which is isomorphic to this expression are presented in figure \ref{f8}.

 \begin{figure}[H]
\centering
\subfloat[Truth table corresponding to the negation of q; that is, $q_1 \land \overline q $.]{
\hspace{1.7in}
\begin{tabular}{c|c||c}
$q$ & $\overline{q}$  &  $q\land \overline{q} $\\
\midrule
0 & 1 & 0   \\  
1 & 0  & 0 \\  
\end{tabular}
\hspace{2in}
\label{f8a}
}
\\
\subfloat[Membership table corresponding to $C$; that is, $C \cap {\xvec{C}}$.]{
\hspace{1.6in}
\begin{tabular}{c|c||c}
$C$ & ${\xvec{C}}$  & $C \cap {\xvec{C}}$ \\
\midrule
0 & 1   & 0\\  
1 & 0   & 0 \\  
\end{tabular}
\hspace{2in}
\label{f8b}}%
\caption{a) Truth table corresponding to $q\land \overline{q} $, and b) membership table corresponding to the set $C \cap {\xvec{C}}$.}
\label{f8}
\end{figure}

The set $C \cap {\xvec{C}}$ is the intersection of the set $C$ and its complement, the set ${\xvec{C}}$. As seen in the membership table corresponding to the set $C \cap {\xvec{C}}$, presented in figure \ref{f8}, no element of the universal set $\mathbb U$ belongs to the set $C \cap {\xvec{C}}$. Therefore, $(C \cap {\xvec{C}})$ is the empty set:  $(C \cap {\xvec{C}} = \varnothing)$. This result is usually regarded as a law in set theory. For this reason, the proposition $q\land \overline{q} $ -- isomorphic to $C \cap {\xvec{C}}$ -- was considered in this section.

\section{The Propositional Calculus of a Canonical\\Fuzz\-y Logic (CFL)}

The CFL mentioned in this article is no different from other fuzzy logics -- including the version of them which is best known and most used, introduced by Lotfi A. Zadeh -- regarding what is understood by ``weight of truth''. The reasons why it is often acceptable and useful to consider that the possible truth values of a proposition are not two -- true and false -- as in BL, have been amply discussed and analyzed critically in specialized literature.

In this article it also will be accepted that within the framework of CFL, each proposition may be assigned a ``weight (or degree) of truth''. Thus, for example, the proposition $q$ may be assigned a weight $w$, which depends on $q$, such that $ 0 \leq w \leq 1$. Reference may be made to the weight of truth using the symbol $w(q)$.

If consideration is given to a set of propositions $q_1, q_2, q_3,\ldots,q_n$, reference may be made to the weight of truth $q_i$, for $i=1,2,3,\ldots,n$, using the symbol $w(q_i)$, and to the weight of truth of $\overline{q}_i$, using the symbol $w(\overline{q}_i)$.

It will be accepted that $w(q_i) + w(\overline{q}_i)=1$; and therefore, $w(q_i)=1-w(\overline{q}_i)$ and $w(\overline{q}_i) = 1-w(q_i)$. The justification for the equality, $w(q_i) + w(\overline{q}_i)=1$, will be provided when considering the law of the excluded middle.

In the applications of CFL, experts in the diverse fields of application should be the ones who determine the weights of truth of the different propositions. Thus, to determine the weight of truth of $q_1$ -- John is tall -- an anthropometric criterion must be used which will depend on the population in which John is considered.  In addition, to specify the weight of truth of $q_2$ -- Ellen is rich -- an economic criterion should be used to establish that quantitative evaluation of the degree of Ellen's wealth, also within the context of a certain reference population.

In figure \ref{f9} and in the explanations given about it in the text below, it will be shown how the weight of truth  $w(q_1 \lor q_2)$ is determined; that is, the weight of truth of the \textit{disjunction} of $q_1$ and $q_2$, once $w(q_1)$ and $w(q_2)$ are known.

\begin{figure}[H]
\centering
\hspace{1in}
\begin{tabular}{c|c|c|c}
$q_1$ & $q_2$ & $q_1\lor q_2$ & $w(q_1 \lor q_2) = \displaystyle\sum_{i=1}^{4}S_i$\\
\midrule
0 & 0 & 0 & $S_1=0$  \\
0 & 1 & 1 & $S_2=w(\overline{q}_1) \cdot w(q_2)$  \\  
1 & 0 & 1 & $S_3=w(q_1) \cdot w(\overline{q}_2)$  \\  
1 & 1  & 1 & $S_4=w(q_1) \cdot w(q_2)$\\  
\end{tabular}
\hspace{1in}
\caption{How to determine $w(q_1 \lor q_2)$. See explanation below.}
\label{f9}

\end{figure}

On the left, or more precisely, in the first three columns of figure \ref{f9}, the truth table for $(q_1 \lor q_2)$ is presented as constructed by using the same procedure as that used in BL propositional calculus. The fourth column of figure \ref{f9} is devoted to the computation of $w(q_1 \lor q_2)$, according to CFL.

The equality $w(q_1 \lor q_2) = \displaystyle\sum_{i=1}^{4}S_i$ specifies that the weight of truth of $q_1 \lor q_2$ -- $w(q_1 \lor q_2)$ -- according to CFL, is obtained by the sum of four terms. In general, there is a term, or addend, corresponding to each row of the truth table considered. In this case, there are 4 rows in the table because the proposition $q_1 \lor q_2$, which is the function of 2 propositions ($q_1$ and $q_2$), is considered.

The first addend, or summand, $S_1$ considered is equal to 0 because there is a 0 in the first row of the column corresponding to $q_1 \lor q_2$. In all cases in which in one row of the column corresponding to the proposition whose weight of truth should be computed there is a 0, the respective addend is equal to 0.

The addends corresponding to the other 3 rows of the column corresponding to $q_1 \lor q_2$ should be computed because in each one of these rows there is a 1, not a 0.

The numerical sequence corresponding to the second row of the truth table for $q_1 \lor q_2$ is as follows: 0, 1, 1. Written out: ``If $q_1$ is false and $q_2$ is true, then $q_1\lor q_2$ is true''. Given that the third digit in the numerical sequence 0, 1, 1 is a 1, the second addend $S_2$ of $w(q_1 \lor q_2)$ should be computed. This addend has 2 factors. The first of them, given that a 0 was assigned to $q_1$, is $w(\overline q_1)$, that is, a truth value of false. The second factor, given that a 1 was assigned to $q_2$, is $w(q_1)$, that is, a truth value of true. Hence, $S_2 = w(\overline q_1)\cdot w(q_2)$.

The numerical sequence corresponding to the third row of the truth table for $q_1 \lor q_2)$ is the following: 1, 0, 1. Written out: ``If $q_1$ is true and $q_2$ is false, then $q_1\lor q_2$ is true''. Given that the third digit in that numerical sequence is a 1, the value of the third addend $S_3$ of $w(q_1 \lor q_2)$ should be computed. This addend has 2 factors. The first of them, given that a 1 was assigned to $q_1$, is $w(q_1)$, that is, a truth value of true. The second factor, given that a 0 was assigned to $q_2$, is $w(\overline{q}_1)$, that is, a truth value of false. Hence, $S_3 = w(q_1)\cdot w(\overline{q}_2)$.

The numerical sequence corresponding to the fourth row of the truth table for $q_1\lor q_2$ is as follows: 1, 1, 1. Written out: ``If $q_1$ is true and $q_2$ is true, then $q_1\lor q_2$ is true''. Given that the third digit in that numerical sequence (1, 1, 1) is a 1, the fourth addend $S_4$ of $w(q_1\lor q_2)$ must be computed. This addend has 2 factors. The first of them, given that a 1 was assigned to $q_1$, is $w(q_1)$, that is, a truth value of true. The second factor, given that a 1 was also assigned to $q_2$, is $w(q_2)$, that is, a truth value of true . Hence, $S_4 = w(q_1)\cdot w(q_1)$.

Therefore, $w(q_1 \lor q_2) = S_1 + S_2 + S_3 + S_4 =
 0 + w(\overline{q}_1) \cdot w(q_2) + w(q_1) \cdot w(\overline{q}_2) + w(q_1) \cdot w(q_2)  = 
 (1-w(q_1)) \cdot w(q_2) + w(q_1) \cdot (1-w(q_2)) + w(q_1) \cdot w(q_2)  =
 w(q_2)-w(q_1)\cdot w(q_2) + w(q_1)-w(q_1)\cdot w(q_2) + w(q_1) \cdot w(q_2)   =
 w(q_1) + w(q_2) - w(q_1)\cdot w(q_2)$.

In conclusion, according to CFL, $w(q_1 \land q_2) = w(q_1) + w(q_2) - w(q_1)\cdot w(q_2)$.

It is important to note that the equality computed according to CFL,\\
$w(q_1 \lor q_2) = w(q_1) + w(q_2) - w(q_1)\cdot w(q_2)$,\\
is different from the equality accepted by the most commonly used fuzzy logics:\\ $w(q_1 \lor q_2) = \max(w(q_1),w(q_2))$.

The same approach is used to obtain the weight of truth $w(q_1 \land q_2)$; that is, the weight of truth according to CFL of the \textit{conjunction} of $q_1$ and $q_2$. Figure \ref{f10} facilitates the understanding of the procedure used here:

\begin{figure}[H]
\centering
\hspace{1in}
\begin{tabular}{c|c|c|c}
$q_1$ & $q_2$ & $q_1\land q_2$ & $w(q_1 \land q_2) = \displaystyle\sum_{i=1}^{4}S_i$\\
\midrule
0 & 0 & 0 & $S_1=0$  \\  
0 & 1 & 0 & $S_2=0$  \\ 
1 & 0 & 0 & $S_3=0$  \\ 
1 & 1  & 1 & $S_4=w(q_1) \cdot w(q_2)$\\  
\end{tabular}
\hspace{1in}
\caption{How to determine $q_1\land q_2$. This is explained below.}
\label{f10}

\end{figure}

Figure \ref{f10} makes it possible to see that $S_4=w(q_1) \cdot w(q_2)$ is the only addend $S_1$, for $i = 1, 2, 3, 4$, in the case of the proposition $q_1\land q_2$ different from 0, according to CFL.

It is important to note that the equality computed according to CFL,
\\ $w(q_1 \land q_2) = w(q_1)\cdot w(q_2)$,\\
is different from the equality accepted by the most commonly used in fuzzy logics:\\
$w(q_1 \land q_2) = \min(w(q_1),w(q_2))$.

In addition, the same approach is used to determine the weights of truth $w(q_1 \to q_2)$, $w(q_2 \to q_1)$ and $w(q_1 \longleftrightarrow q_2)$, as specified in figure \ref{f11} .

 \begin{figure}[H]
\centering
\subfloat[How to determine $w(q_1 \to q_2)$.]{
\hspace{1in}
\begin{tabular}{c|c|c|c}
$q_1$ &  $q_2$ & $q_1 \to q_2$  & $w(q_1 \to q_2) =  \displaystyle\sum_{i=1}^{4}S_i$\\
\midrule
0 & 0 & 1 & $S_1=w(\overline{q}_1) \cdot w(\overline{q}_2)$  \\  
0 & 1 & 1 & $S_2=w(\overline{q}_1) \cdot w(q_2)$  \\ 
1 & 0 & 0 & $S_3=0$  \\  
1 & 1  & 1 & $S_4=w(q_1) \cdot w(q_2)$\\ 
 
\end{tabular}
\hspace{1in}
\label{f11a}
}
\\
\subfloat[How to determine $w(q_2 \to q_1)$.]{
\hspace{1in}
\begin{tabular}{c|c|c|c}
$q_1$ &  $q_2$ & $q_2 \to q_1$  & $w(q_2 \to q_1) = \displaystyle\sum_{i=1}^{4}S_i$  \\
\midrule
0 & 0 & 1 & $S_1=w(\overline{q}_1) \cdot w(\overline{q}_2)$  \\ 
0 & 1 & 0 & $S_2=0$  \\ 
1 & 0 & 1 & $S_3=w(q_1) \cdot w(\overline{q}_1)$  \\  
1 & 1  & 1 & $S_4=w(q_1) \cdot w(q_2)$\\ 
 
\end{tabular}
\hspace{1in}
\label{f11b}}%
\\
\subfloat[How to determine $w(q_1 \longleftrightarrow q_2)$.]{
\hspace{1in}
\begin{tabular}{c|c|c|c}
$q_1$ &  $q_2$ & $q_1 \longleftrightarrow q_2$  & $w(q_1 \longleftrightarrow q_2) = \displaystyle\sum_{i=1}^{4}S_i$  \\
\midrule
0 & 0 & 1 & $S_1=w(\overline{q}_1) \cdot w(\overline{q}_2)$  \\  
0 & 1 & 0 & $S_2=0$  \\
1 & 0 & 0 & $S_3=0$  \\ 
1 & 1  & 1 & $S_4=w(q_1) \cdot w(q_2)$\\ 
 
\end{tabular}
\hspace{1in}
\label{f11c}}%
\caption{How to determine a) $w(q_1 \to q_2)$, b) $w(q_2 \to q_1)$, and c) $w(q_1 \longleftrightarrow q_2)$.}
\label{f11}
\end{figure}

Recall that, for any $q_i$, for $i=1,2,3,\ldots$, $w(q_i) + w(\overline{q}_i)=1$. Therefore, according to the information provided in figure \ref{f11}, the following results are obtained:\\
\\
$w(q_1 \to q_2)=w(\overline{q}_1) \cdot w(\overline{q}_2) +  w(\overline{q}_1) \cdot w(q_2) + 0 + w(q_1) \cdot w(q_2) = \\ w(\overline{q}_1)\cdot (w(q_2)+w(\overline{q}_2)) + w(q_1)\cdot w(q_2) = w(\overline{q}_1) + w(q_1)\cdot w(q_2)  = \\ 1-w(q_1)+w(q_1)\cdot w(q_2)$\\ 
\\
$w(q_2 \to q_1)=w(\overline{q}_1) \cdot w(\overline{q}_2) + 0 + w(q_1) \cdot w(\overline{q}_1) + w(q_1) \cdot w(q_2) = \\ w(\overline{q}_2)\cdot w(\overline{q}_1 + w(q_1)) + w(q_1)\cdot w(q_2) = \\ 1-w(q_2)+w(q_1)\cdot w(q_2)$\\ 
\\
$w(q_1 \longleftrightarrow q_2)=w(\overline{q}_1) \cdot w(\overline{q}_2)  + 0 + 0 + w(q_1) \cdot w(q_2) = \\ (1 - w(q_1))\cdot(1-w(q_2)) +  w(q_1)\cdot w(q_2) = \\ 1- w(q_2) - w(q_1) + w(q_1)\cdot w(q_2) + w(q_1)\cdot w(q_2) = \\  1-w(q_1) - w(q_2)+ 2\cdot w(q_1)\cdot w(q_2) $.\\

Hence, according to CFL,\\ 
$w(q_1 \to q_2) = 1-w(q_1)+w(q_1)\cdot w(q_2) \\ w(q_2 \to q_1)= 1-w(q_2) + w(q_1)\cdot w(q_2) \ \\ w(q_1 \longleftrightarrow q_2)= 1-w(q_1) - w(q_2)+ 2\cdot w(q_1)\cdot w(q_2)  $.

Consideration is given below to the law of the excluded middle according to the propositional calculus of CFL.

\begin{figure}[H]
\centering
\hspace{1in}
\begin{tabular}{c|c|c|c}
$q$ & $\overline{q}$ & $q\lor \overline{q}$ & $w(q \lor \overline{q}) = \displaystyle\sum_{i=1}^{2}S_i$\\
\midrule
0 & 1& 1 & $S_1=w(\overline{q})$  \\  
1 & 0  & 1 & $S_2=w(q) $\\  
\end{tabular}
\hspace{1in}
\caption{How to determine $w(q \lor \overline{q})$.}
\label{f12}
\end{figure}

According to figure \ref{f12}, the following result is obtained:\\
$w(q\lor \overline{q}) = S_1 + S_2 = w(\overline{q}) + w(q) =1 - w(q) + w(q) = 1$
\\
In this case one might wonder about how $S_1$ and $S_2$ were specified. The justification is as follows: $q\lor \overline{q}$ is a propositional function of a sole proposition $q$, and in general, for any propositional function of $n$ propositions, for $n=1,2,3\ldots$, each $S_i$ is equal to the product of $n$ factors, each of which is a certain weight of truth.

In figure \ref{f13} indications are given on how to determine the weight of truth of the proposition $q\land \overline{q}$ according to the propositional calculus of CFL.

\begin{figure}[H]
\centering
\hspace{1in}
\begin{tabular}{c|c|c|c}
$q$ & $\overline{q}$ & $q\land \overline{q}$ & $w(q \land \overline{q}) = \displaystyle\sum_{i=1}^{2}S_i$\\
\midrule
0 & 1 & 0 & $S_1=0 $  \\  
1 & 0  & 0 & $S_2=0  $\\  
\end{tabular}
\hspace{1in}
\caption{How to determine $w(q \land \overline{q})$.}
\label{f13}
\end{figure}

According to figure \ref{f13}, the following result is obtained:\\
$w(q\land \overline{q}) = S_1 + S_2 = 0$.

The reason why $S_1 = 0$ here is as follows: In the first row of the column corresponding to $q\land \overline{q}$, there is a 0. The reason why $S_2 = 0$ here is as follows: In the second row of the column corresponding to $q\land \overline{q}$, there also is a 0. Given that the weight of truth of $q\land \overline{q}$ (i. e., $w(q\land \overline{q})$) is equal to 0, regardless of the weight of truth of $q$ (i. e., $w(q)$), it is considered that $q\land \overline{q}$ is a contradiction of propositional calculus of CFL.

Another tautology, one of De Morgan's laws, will be considered below according to the propositional calculus of CFL.

\begin{figure}[H]
\centering
\hspace{1in}
\begin{tabular}{c|c|c|c}
$q_1$ & $q_2$ & $ (\overline{q_1 \lor q_2}) \longleftrightarrow (\overline{q}_1 \land \overline{q}_2)$ & $w((\overline{q_1 \lor q_2}) \longleftrightarrow (\overline{q}_1 \land \overline{q}_2)) = \displaystyle\sum_{i=1}^{4}S_i$\\
\midrule
0 & 0 & 1 & $S_1=w(\overline{q}_1) \cdot w(\overline{q}_2)$  \\  
0 & 1 & 1 & $S_2=w(\overline{q}_1) \cdot w(q_2) $  \\ 
1 & 0 & 1 & $S_3=w(q_1) \cdot w(\overline{q}_1)$  \\ 
1 & 1  & 1 & $S_4=w(q_1) \cdot w(q_2)$\\
\end{tabular}
\hspace{1in}
\caption{How to determine $w((\overline{q_1 \lor q_2}) \longleftrightarrow (\overline{q}_1 \land \overline{q}_2))$.}
\label{f14}
\end{figure}

According to figure \ref{f14},  the following result is obtained:\\
$w((\overline{q_1 \lor q_2}) \longleftrightarrow (\overline{q}_1 \land \overline{q}_2)) = w(\overline{q}_1) \cdot w(\overline{q}_2)+w(\overline{q}_1) \cdot w(q_2)+w(q_1) \cdot w(\overline{q}_2)+w(q_1) \cdot w(q_2) = w(\overline{q}_1)\cdot(w(\overline{q}_2)+w(q_2)) +  w(q_1)\cdot(w(\overline{q}_2)+w(q_2)) = w(\overline{q}_1) +w(q_1) = 1$.

Recall that a tautology of BL propositional calculus is a propositional function of $n$ propositions, for $n=1,2,3\ldots$, which is true regardless of the truth value -- true or false -- of one of the $n$ propositions.

In general, any tautology of BL propositional calculus is also a tautology -- or law -- of  CFL propositional calculus, in the sense that its weight of truth is equal to 1, regardless of the weight of truth of each of the $n$ propositions of which it is a propositional function.

A negation of any tautology of BL propositional calculus is a contradiction. Recall that in that calculus a contradiction is any propositional function of $n$ propositions, for $n=1,2,3\ldots$, which is false regardless of the truth value -- true or false -- of each of those $n$ propositions.

In general, any contradiction of propositional calculus of BL is also a contradiction of propositional calculus of CFL, in the sense that its weight of truth is equal to 0, regardless of the weight of truth of each of the $n$ propositions of which it is a propositional function.

As a consequence of the characterizations given of tautologies and contradictions, both of propositional calculus of BL and of propositional calculus of CFL, it is obvious that in both types of calculus, a) the negation of any tautology is a contradiction, and b) the negation of any contradiction is a tautology.

All the results that can be obtained with the propositional calculus of BL can also be obtained with the propositional calculus of CFL. Indeed, in the former, each proposition can have only one of two truth values: true or false. If each proposition of this calculus is reinterpreted within the framework of propositional calculus of CFL, restricting the weights of truth to only two of them (the weight of truth equal to 1, or the weight of truth equal to 0), then the second propositional calculus -- that of CFL -- makes it possible to obtain all the results achieved by using the first of these calculi -- that of BL. Thus it can be concluded that the first type of propositional calculus is one particular case, which could be viewed as a ``limit case'' of the second type of calculus given that 0 and 1 are the minimum and maximum numerical values respectively, which in the second calculation can be assigned to the weight of any proposition.

\section{Fuzzy Sets According to CFL}

In classical set theory, one possible way of characterizing a finite set is to specify each element belonging to it. Thus, the set $C_1$ to which the natural numbers 7, 17 and 528 belong can be characterized as follows:

\begin{equation}
C_1 = \{7, 17, 528\} 
\end{equation}

In some applications of this theory which mention different extra-mathematical entities, such as human beings, an expression such as the following does not suffice to characterize a set unambiguously:
 
\begin{equation}
C_2 = \{x_1, x_2, x_3\}    
\end{equation}

In effect, the characterization of $C_2$ must also specify which people are referred to by the symbols $x_1, x_2, x_3$. However, difficulties may appear in the applications to different fields of knowledge. Some of these difficulties have led to the introduction of logics different from BL. Suppose that $C_2$ is the ``set of tall men in a certain population''. In the case considered in (1), there is no doubt as to whether each element of the set $N$ of the natural numbers belongs, or does not belong, to the set $C_1$. Each natural number either belongs or does not belong to $C_1$. Therefore, for example, given the characterization of $C_1$, it is known that the natural number 287 does not belong to $C_1$ and that the natural number 528 does belong to $C_1$. On the other hand, in case (2) it could be useful to introduce the notion of degree or ``weight of membership'' in the set $C_2$ for each element considered. This was covered in section 4 above. Considerations of this type were precisely those that led to the development of fuzzy set theory.

Just as in classical set theory, in CFL fuzzy set theory, when applying the latter to the treatment of a given topic, first the characterization of the universal set $\mathbb U$ to which all the elements belong is introduced, and one must specify the elements to which reference is made.  It will be accepted that each element belongs to a $\mathbb U$ with a weight of membership equal to 1. In this article only finite fuzzy universal sets will be discussed. The notation to be used is adequate for this restriction. As this will be covered in other articles, that limitation can be eliminated to make it possible to operate, using the CFL approach, with both countable and uncountable infinite universal sets.

In classical set theory, each element belonging to $\mathbb U$ does not belong to any set $C$ different from $\mathbb U$. In effect, if each element of those belonging to $\mathbb U$ belonged to $C$, $C$ would not be different from $\mathbb U$, but rather equal to it; that is, the equality $C = \mathbb U$ would be valid.

In fuzzy set theory, according to CFL, each cf elements belonging to $\mathbb U$ does belong to each set characterized within the framework of a given $\mathbb U$. However, if that set is different from $\mathbb U$, then at least one of those elements belongs to it with a weight different from 1. As justified below, each element belonging to $\mathbb U$ belongs to the complement of $\mathbb U$ ($\xvec{\mathbb U}$), which is equal to the empty set $\varnothing$, with a weight of membership equal to 0.

Consider any element $x_i$ whatsoever, for $i = 1,2,3,\dots,n$, belonging to $\mathbb U$. That membership may be expressed as $x_i \in \mathbb U$. (Recall that the symbol $\in$ is used to indicate the membership of a given element in a set.)

The weight of membership of $x_i$ in $\mathbb U$ is symbolized as $w(x_i \in \mathbb U)$. In this case, that weight of membership is equal to 1: $w(x_i \in \mathbb U)=1$.

In general, for the weight of membership of any element $x_i$ of $\mathbb U$ in a set $C_j$, for $j = 1,2,3,\dots$, characterized within the framework of that $\mathbb U$, it will be accepted that its numerical value is between 0 and 1: $0 \leq w(x_i \in C_j) \leq 1$.

Notice that ($x_i \in \mathbb U$) and ($x_i \in C_j$) are, respectively, the propositions ``The element $x_i$ belongs to $\mathbb U$'' and ``The element $x_i$ belongs to $C_j$''. Therefore, the weight of membership of any element in a set may be interpreted also in CFL set theory as the weight of truth of the proposition which affirms the membership of that element in the set considered.

The propositions $x_i \in \mathbb U$ and $x_i \in \xvec{\mathbb U}$  will be expressed, respectively, in an abbreviated way, as $q_{i,\mathbb U}$ and $q_{i,\xvecsub{\mathbb U}}$. The first subscript in both cases specifies which element is mentioned and the second subscript specifies which sets -- $\mathbb U$ and $ \xvec{\mathbb U}$, respectively -- the element is said to belong to. Given that $ \xvec{\mathbb U} = \varnothing$, $q_{i,\xvecsub{\mathbb U}} = q_{i,\varnothing}$.

The propositions ($x_i \in C_j$)  and ($x_i \in \xvec{C}_j$) will be expressed, respectively, in an abbreviated way, as $q_{i,j}$, and $q_{i,\xvecsub{j}}$. The first subscript in both cases specifies unambiguously which sets $C_j$ and $\xvec{C}_j$, respectively, the element belongs to.

As seen above, each proposition of this type has a weight of truth. Consider, for example, the proposition $q_{i,j}$. Its weight of truth is symbolized as $w(q_{i,j})$. That proposition $q_{i,j}$ will be regarded as a proposition characteristic of CFL propositional calculus. Thus, $w(q_{i,j}) +  w(\overline{q}_{i,j}) = 1$, and
\begin{equation}
w(\overline{q}_{i,j}) = 1 - w(q_{i,j}).               
\end{equation}

It will be accepted that the following equality is valid: $w(q_{i,\xvecsub{j}}) = w(\overline{q}_{i,j})$. Therefore, due to (3), 
\begin{equation}
w(q_{i,\xvecsub{j}}) = 1 - w(q_{i,j}).         
\end{equation}
Written out, equation (4) can be expressed as follows: The weight of truth of $q_{i,\xvecsub{j}}$ is equal to the difference between 1 and the weight of truth of $q_{i,j}$, or in an equivalent way, the weight of membership of any element belonging to the complement of any set $C_j$ (that is, $\xvec{C}_j$) is equal to 1 minus the weight of membership of that element in $C_j$. If this result is applied to the universal set $\mathbb U$, the following is obtained: $w(q_{i,\xvecsub{\mathbb U}}) = w(q_{i,\varnothing}) = 1-w(q_{i,\mathbb U}) = 1-1 = 0 $. As previously stated, $w(q_{i,\varnothing}) = 0$.

If desired, to recall that one is operating with fuzzy sets, according to CFL, the superscript $f$, for example, for \textit{f}uzzy may be added on the left to each set of this type, in the following way:\\
$~^f\mathbb U$; $~^f\varnothing$; $~^fC_5$; $~^f\xvec{C}_9$.\\
That will not be done here. In many cases it is unnecessary because the type of set with which one is operating is clear in context or by clarification.

In applied logic one often selects a universal set $\mathbb U$ such that all the elements belonging to it are the same type. In this section, to provide examples, according to CFL, of some operations with sets, attention  will be given below to a  $\mathbb U$ such that each of the 5 elements belonging to it is a person -- a human being. Consideration will also be given to 2 of the subsets of that $\mathbb U$: $C_1$ and $C_2$.

$C_1$ and $C_2$ will be characterized as follows:\\
$C_1$: set of chess players, and\\
$C_2$ set of wealthy people.

$\mathbb U$, to which it will be admitted that persons $x_1$, $x_2$, $x_3$, $x_4$ and $x_5$ belong, was selected in a way such that $x_2$ is a distinguished professional chess player, a grandmaster, with the possibility of becoming the next world champion in chess. According to experts in chess, it may be considered that $w(q_{2,1}) = 1$.

The element $x_5$ is a chess fan whose level, according to experts in chess, is that of a third class player: $w(q_{5,1}) = 0.4$.

Regarding $x_1$, $x_3$ and $x_4$, none of these people know anything about chess rules. They do not know, for instance, how to move each one of the chess pieces, according to the rules of the game. According to chess experts, $w(q_{1,1}) = w(q_{3,1}) = w(q_{4,1}) = 0$.

In another area, let us admit that experts in economics and finances came to the following conclusions:\\
$w(q_{1,2}) = 0.9$; $w(q_{2,2}) = 0.8$; $w(q_{3,2}) = 0.7$; $w(q_{4,2}) = 0$; and $w(q_{5,2}) = 0.6$.

For clarity and to facilitate the comprehension of the operations to be carried out with the objective of providing some examples of them, the different fuzzy sets to be discussed will be expressed below as column vectors. Each row of these column vectors begins with the symbol $x_i$ for $i = 1,2,3,4,5$, corresponding to the element considered; and to the right of that symbol and separated from it by a semi-colon (;), the equation which establishes the corresponding weight of truth is provided.

\begin{equation*}
\mathbb U = \begin{array}{|cc|}
x_1;  & w(q_{1,\mathbb U}) =1   \\
x_2;  & w(q_{2,\mathbb U})=1 \\
x_3;  & w(q_{3,\mathbb U}) =1 \\
x_4;  & w(q_{4,\mathbb U}) =1 \\
x_5;  & w(q_{5,\mathbb U}) =1 \\
\end{array}\quad;\qquad
\xvec{\mathbb U} =\varnothing = \begin{array}{|cc|}
x_1;  & w(q_{1,\varnothing}) =0   \\
x_2;  & w(q_{2,\varnothing})=0 \\
x_3;  & w(q_{3,\varnothing}) =0 \\
x_4;  & w(q_{4,\varnothing}) =0 \\
x_5;  & w(q_{5,\varnothing}) =0 \\
\end{array}\end{equation*}

\begin{equation*}
C_1 = \begin{array}{|cc|}
x_1;  & w(q_{1,1}) =0   \\
x_2;  & w(q_{2,1})=1 \\
x_3;  & w(q_{3,1}) =0 \\
x_4;  & w(q_{4,1}) =0 \\
x_5;  & w(q_{5,1}) =0.4 \\
\end{array}\quad;\qquad
\xvec{C}_1 = \begin{array}{|cc|}
x_1;  & w(q_{1,\xvecsub{1}}) =1   \\
x_2;  & w(q_{2,\xvecsub{1}})=0 \\
x_3;  & w(q_{3,\xvecsub{1}}) =1 \\
x_4;  & w(q_{4,\xvecsub{1}}) =1 \\
x_5;  & w(q_{5,\xvecsub{1}}) =0.6 \\
\end{array}\end{equation*}
\\
For each element $x_i$, for $i = 1,2,3,4,5$, the numerical values of $w(q_{i,\xvecsub{1}})$ were obtained from $w(q_{i,1})$ by using (4).

\begin{equation*}
C_2 = \begin{array}{|cc|}
x_1;  & w(q_{1,2}) =0.9   \\
x_2;  & w(q_{2,2})=0.8 \\
x_3;  & w(q_{3,2}) =0.7 \\
x_4;  & w(q_{4,2}) =0 \\
x_5;  & w(q_{5,2}) =0.6 \\
\end{array}\quad;\qquad
\xvec{C}_2 = \begin{array}{|cc|}
x_1;  & w(q_{1,\xvecsub{2}}) =0.1  \\
x_2;  & w(q_{2,\xvecsub{2}})=0.2 \\
x_3;  & w(q_{3,\xvecsub{2}}) =0.3 \\
x_4;  & w(q_{4,\xvecsub{2}}) =1 \\
x_5;  & w(q_{5,\xvecsub{2}}) =0.6 \\
\end{array}\end{equation*}
\\
Again, (4) was used to obtain numerical values of $w(q_{i,\xvecsub{2}})$, from the numerical valus of $w(q_{i,2})$.

Some operations carried out with the fuzzy sets specified are presented below to emphasize that all the laws of classical set theory are conserved -- that is, also valid -- in CFL fuzzy set theory.

In classical set theory, for any set $C$, the following law is valid: $C \cup \xvec{C} = \mathbb U$. Written out: The union of any set $C$ whatsoever and its complement $\xvec{C}$ is equal to the universal set $\mathbb U$.

Consider any universal set $\mathbb U$ according to CFL set theory. It will be proven that if the operation of union of that universal set and its complement $\xvec{\mathbb U}$ is carried out, the following result is obtained:\\
$\mathbb U \cup \xvec{\mathbb U} = \mathbb U \cup \varnothing = \mathbb U$

In effect, the weight of truth $w(q_{i,\mathbb U} \lor q_{i,\varnothing})$ can be computed for each $x_i$, for $i = 1,2,3,\ldots, n$, belonging to that $\mathbb U$:\\
$w(q_{i,\mathbb U} \lor q_{i,\varnothing}) = w(q_{i,\mathbb U} \lor q_{i,\xvecsub{\mathbb U}}) = w(q_{i,\mathbb U} \lor \overline{q}_{i,\mathbb U}) = w(q_{i,\mathbb U})+w(\overline{q}_{i,\mathbb U})  = w(q_{i,\mathbb U}) + 1- w(q_{i,\mathbb U}) = 1$.

Given that for any element $x_i$, for $i = 1,2,3,\ldots, n$, belonging to $\mathbb U \cup \varnothing$, the corresponding weight of truth $w(q_{i,\mathbb U} \lor q_{i,\varnothing})$ is equal to 1, the union of $\mathbb U$ and $\varnothing$ ($\mathbb U \cup \varnothing$) is equal to the set $\mathbb U$: $\mathbb U \cup \varnothing = \mathbb U$.

Consider any subset $C_j$, for $j = 1,2,3,\ldots$, of a $\mathbb U$, according to CFL. In general, for any element $x_i$, for $i = 1,2,3,\ldots$, $n$, belonging to ($C_j \cup \xvec{C_j}$), the following equality is valid:\\
$w(q_{i,j} \lor q_{i,\xvecsub{j}}) = w(q_{i,j} \lor \overline{q}_{i,j}) = w(q_{i,j}) + w(\overline{q}_{i,j}) =1$.\\
Therefore, ($C_j \cup \xvec{C}_j = \mathbb U$).

In particular, it can be easily verified that for the set $\mathbb U$ of the 5 people mentioned, from which the subset $C_1$ of chess players and  the subset $C_2$ of wealthy people were considered, the following equalities are valid:\\
\\
$\mathbb U \cup \xvec{\mathbb U} = \mathbb U \cup \varnothing = \mathbb U$\\
\\
$C_1 \cup \xvec{C}_1 = \mathbb U$\\
\\
$C_2 \cup \xvec{C}_2 = \mathbb U$\\
 
 For any element $x_i$ of any $\mathbb U$, according to CFL, the following equalities are valid:\\
 $w(q_{i,\mathbb U} \land q_{i,\varnothing}) = w(q_{i,\mathbb U} \land  q_{i,\xvecsub{\mathbb U}}) =  w(q_{i,\mathbb U} \land \overline{q}_{i,\mathbb U})  = 0$.\\
In effect, according to what was discussed in section 4, the weight of truth of the conjunction of a proposition and the negation of that proposition is equal to 0. Therefore,\\
$\mathbb U \cap \varnothing=\varnothing$.

Reasoning of this same type is valid for any element $x_i$ of any subset $C_j$ of that $\mathbb U$: $w(q_{i,j} \land q_{i,\xvecsub{j}}) = w(q_{i,j} \land \overline{q}_{i,j}) = 0$.
 
In particular, for the set $\mathbb U$ of 5 people of which $C_1$ and $C_2$ are subsets, it can easily be verified that the following equalities are valid:\\
\\
$\mathbb U \cap \xvec{\mathbb U} = \mathbb U \cap \varnothing = \varnothing$\\
\\
$C_1 \cap \xvec{C}_1 = \varnothing$\\
\\
$C_2 \cap \xvec{C}_2 = \varnothing$\\

Consideration will be given below to the union of $C_1$ and $C_2$ ($C_1 \cup C_2$) and to the intersection of those sets ($C_1 \cap C_2$).

\begin{align*}
(C_1 \cup C_2) = \begin{array}{|cc|}
x_1;  & w(q_{1,{1}}) =0   \\
x_2;  & w(q_{2,{1}})=1 \\
x_3;  & w(q_{3, {1}}) =0 \\
x_4;  & w(q_{4, {1}}) =0 \\
x_5;  & w(q_{5, {1}}) =0.4 \\
\end{array}  \hspace{.1in} \cup \hspace{.1in} \begin{array}{|cc|}
x_1;  & w(q_{1,2}) =0.9   \\
x_2;  & w(q_{2,2})=0.8 \\
x_3;  & w(q_{3,2}) =0.7 \\
x_4;  & w(q_{4,2}) =0 \\
x_5;  & w(q_{5,2}) =0.6 \\
\end{array} =\\ \begin{array}{|cc|}
x_1;  & w(q_{1,1} \lor q_{1,{2}}) =0.9  \\
x_2;  & w(q_{2,1} \lor q_{2,{2}})= 1 \\
x_3;  & w(q_{3,1} \lor q_{3,{2}}) =0.7 \\
x_4;  & w(q_{4,1} \lor q_{4,{2}}) =0 \\
x_5;  & w(q_{5,1} \lor q_{5,{2}}) =0.76 \\
\end{array}\end{align*}
\\
The five weights of truth for $w(q_{i,1} \lor q_{i,2})$, for $i = 1, 2, 3, 4,5$, were computed as specified in section 4.

\begin{align*}
(C_1 \cap C_2) = \begin{array}{|cc|}
x_1;  & w(q_{1,{1}}) =0   \\
x_2;  & w(q_{2,{1}})=1 \\
x_3;  & w(q_{3, {1}}) =0 \\
x_4;  & w(q_{4, {1}}) =0 \\
x_5;  & w(q_{5, {1}}) =0.4 \\
\end{array} \hspace{.1in}  \cap \hspace{.1in}  \begin{array}{|cc|}
x_1;  & w(q_{1,2}) =0.9   \\
x_2;  & w(q_{2,2})=0.8 \\
x_3;  & w(q_{3,2}) =0.7 \\
x_4;  & w(q_{4,2}) =0 \\
x_5;  & w(q_{5,2}) =0.6 \\
\end{array} = \\  \begin{array}{|cc|}
x_1;  & w(q_{1,1} \land q_{1,{2}}) =0  \\
x_2;  & w(q_{2,1} \land q_{2,{2}})= 0.8 \\
x_3;  & w(q_{3,1} \land q_{3,{2}}) =0 \\
x_4;  & w(q_{4,1} \land q_{4,{2}}) =0 \\
x_5;  & w(q_{5,1} \land q_{5,{2}}) =0.24 \\
\end{array}\end{align*}

The five weights of $w(q_{i,1} \land q_{i,2})$, for $i = 1, 2, 3, 4,5$, were computed as specified in section 4.

Consider two subsets whatsoever $C_j$ and $C_k$ of any fuzzy universal set $\mathbb U$, according to CFL. The set $C_l$ is characterized as:\\
$C_l=( (\xvec{C_j \cup  C_k})  \naturaltolr (\xvec{C}_j \cap \xvec{C}_k) )$.

If $C_l$ is equal to $\mathbb U$ (that is, $C_l=\mathbb U$) in the sense that any element $x_i$ of $\mathbb U$ belongs to $C_l$ with a weight of membership equal to 1 (or in an equivalent way, such that the weight of truth $w(q_{i,l})$ is equal to 1, that is, $w(q_{i,l})=1$), then it is proven that one of De Morgan's laws is valid also for the CFL fuzzy sets.

The membership table of any element of $\mathbb U$ belonging to $C_l$ is shown in figure 15.

\begin{figure}[H]
\centering
\hspace{1in}
\begin{tabular}{c|c|c|c}
$C_j$ & $C_k$ & $C_l = ( (\xvec{C_j \cup  C_k})  \naturaltolr (\xvec{C}_j \cap \xvec{C}_k) )   $ & $w(q_{i,l}) = \displaystyle\sum_{i=1}^{4}S_i$\\
\midrule
0 & 0 & 1 & $S_1=w(\overline{q}_{i,j}) \cdot w(\overline{q}_{i,k})$  \\  
0 & 1 & 1 & $S_2=w(\overline{q}_{i,j}) \cdot w(q_{i,k}) $  \\ 
1 & 0 & 1 & $S_3=w(q_{i,j}) \cdot w(\overline{q}_{i,k})$  \\ 
1 & 1  & 1 & $S_4=w(q_{i,j}) \cdot w(q_{i,k})$\\
\end{tabular}
\hspace{1in}
\caption{Table of membership of any element $x_i$ of $\mathbb U$ belonging to $C_l$. As shown above, it may be interpreted also as a particular truth table.}
\label{f15}
\end{figure}

As shown in figure 15,\\
$w(q_{i,l})= \displaystyle\sum_{i=1}^{4}S_i =  w(\overline{q}_{i,j}) \cdot w(\overline{q}_{i,k}) +  w(\overline{q}_{i,j}) \cdot w(q_{i,k})  + w(q_{i,j}) \cdot w(\overline{q}_{i,k}) + w(q_{i,j}) \cdot w(q_{i,k}) = w(\overline{q}_{i,j}) \cdot (w(\overline{q}_{i,k}) + w(q_{i,k})) + w({q}_{i,j}) \cdot (w(\overline{q}_{i,k}) + w(q_{i,k}))   = w(\overline{q}_{i,j})\cdot 1 + w(q_{i,j}) \cdot 1 = w(\overline{q}_{i,j})+  w(q_{i,j}) = 1$.

Therefore, $C_l = \mathbb U$; De Morgan's law mentioned above is also valid for CFL fuzzy sets.

In a future article, consideration will be given to the following result: Any law of classic set theory is also valid in CFL fuzzy set theory.

If for each fuzzy set characterized within the framework of a specific $\mathbb U$, the numerical value of the weight of membership in that set is known or may be computed for each element of $\mathbb U$, which may be reinterpreted as the weight of truth of the proposition which affirms that membership, then a convention that simplifies the characterization of these sets can be accepted. In effect, to characterize any of these sets, the following dyad is provided for each of six elements: 1) the denomination of the element considered, and 2) the numerical value of its weight of membership in that set. This convention will be used below to emphasize an important one-to-one correspondence between each ``classical'' set and a particular fuzzy set, according to CFL. This correspondence is obtained in the following way: Each element of $\mathbb U$ belonging to the ``classical'' set considered will belong also to the corresponding fuzzy set with a weight of membership of 1 and each element of $\mathbb U$ not belonging to the ``classical'' set will belong to the corresponding fuzzy set with a weight of membership equal to 0. 

In the following sequence of correspondences, in which the ``classical'' sets were represented to the left and the fuzzy sets to the right, the convention specified has been used.

\begin{equation*}
\mathbb  U =  \{1,2,3,4,5,6,7\} \quad \text{corresponds to} \quad \mathbb U = \begin{array}{|cc|}
1;  &1   \\
2;  &  1 \\
3;  & 1 \\
4;  & 1 \\
5;  & 1 \\
6;  & 1 \\
7;  & 1 \\
\end{array}   
\end{equation*}
\\
\begin{equation*}
\xvec{\mathbb  U}  =\varnothing =  \{\quad\}  \quad \text{corresponds to} \quad \xvec{\mathbb  U}  =\varnothing  = \begin{array}{|cc|}
1;  & 0   \\
2;  &  0 \\
3;  & 0\\
4;  & 0 \\
5;  & 0 \\
6;  & 0 \\
7;  & 0 \\
\end{array}   
\end{equation*}
\\
\begin{equation*}
C_1 =  \{2,3,5,7\}  \quad \text{corresponds to} \quad C_1 = \begin{array}{|cc|}
1;  & 0  \\
2;  &  1 \\
3;  & 1 \\
4;  & 0 \\
5;  & 1 \\
6;  & 0 \\
7;  & 1 \\
\end{array}   
\end{equation*}
\\
\begin{equation*}
\xvec{C}_1 =  \{1,4,6\}  \quad \text{corresponds to} \quad \xvec{C}_1 = \begin{array}{|cc|}
1;  &1   \\
2;  &  0 \\
3;  & 0 \\
4;  & 1 \\
5;  & 0 \\
6;  & 1 \\
7;  & 0 \\
\end{array}   
\end{equation*}
\\
\begin{equation*}
C_2 =  \{1,3,4,7\}  \quad \text{corresponds to} \quad C_2 = \begin{array}{|cc|}
1;  & 1  \\
2;  &  0 \\
3;  & 1 \\
4;  & 1 \\
5;  & 0 \\
6;  & 0 \\
7;  & 1 \\
\end{array}   
\end{equation*}
\\
\begin{equation*}
\xvec{C}_2 =  \{2,5,6\}  \quad \text{corresponds to} \quad \xvec{C}_2 = \begin{array}{|cc|}
1;  &0  \\
2;  &  1 \\
3;  & 0 \\
4;  & 0 \\
5;  & 1 \\
6;  & 1 \\
7;  & 0 \\
\end{array}   
\end{equation*}
\\
\begin{equation*}
(C_1 \cup C_2)  =  \{1,2,3,4,5,7\}  \quad \text{corresponds to} \quad (C_1 \cup C_2) = \begin{array}{|cc|}
1;  & 1  \\
2;  &  1 \\
3;  & 1 \\
4;  & 1 \\
5;  & 1 \\
6;  & 0 \\
7;  & 1 \\
\end{array}   
\end{equation*}
\\
\begin{equation*}
(C_1 \cap C_2)  =  \{3,7\}  \quad \text{corresponds to} \quad (C_1 \cap C_2) = \begin{array}{|cc|}
1;  & 0  \\
2;  &  0 \\
3;  & 1 \\
4;  & 0 \\
5;  & 0 \\
6;  & 0 \\
7;  & 1 \\
\end{array}   
\end{equation*}\\

\section{Discussion and Perspectives}

Some basic results have been presented regarding a variant of fuzzy logic: canonical fuzzy logic (CFL).

Two of the main objectives of this article, and of others to be published in the future, are:

1. to show the naturally existing continuity for each logic discussed (that of BL and diverse non-classical logics) between calculi of different levels, i.e., between propositional calculus and predicate calculus, which can be presented using set theory terminology; and

2. to specify how BL can be considered a particular ``limit case'' of the other logics to be considered. Thus, for example, BL propositional calculus can be regarded as a particular limit case of CFL propositional logic if in the latter the weight of truth of any proposition considered is restricted to only 2 possible numerical values, 0 or 1. Likewise, as will be discussed in future articles, classical set theory is isomorphic to a particular case of CFL set theory, that in which the weight of truth $w(x \in  C)$ can have only one of two numerical values: 0 or 1. That case refers to the weight of truth of the proposition which indicates the membership of any element $x$ of the universal set $\mathbb U$ considered, belonging to any set $C$ characterized within the framework of that $\mathbb U$.

As seen above, the treatment given in this article to CFL set theory is still quite incomplete. For example, important operators which would extend its scope and provide versatility to its applications remain to be introduced. In addition, it is necessary to specify how CFL operates with fuzzy sets characterized within the frameworks corresponding to infinite universal fuzzy sets.\\

\end{document}